\documentclass[12pt]{article}
\usepackage{amssymb}

\usepackage{amsfonts,amssymb,eucal,amsmath}
\title{
A New Extension of Serrin's Lower Semicontinuity Theorem }
\author{
{Hu Xiaohong$^{1,2}$ } and \quad \quad {Zhang Shiqing$^{1}$}
\footnote{Corresponding author: huxh@cqupt.edu.cn and
zhangshiqing@msn.com}
\thanks{Supported partially by NSF of China. }
\\
{\small $^{1}$Department of Mathematics, Sichuan University,}\\
{\small Chengdu 610064, PR China}\\
{\small $^{2}$ Department of Mathematics and Physics,}\\
{\small Chongqing University of
Posts and Telecommunications,}\\
 {\small Chongqing 400065, PR China}\\
 }
\date{}
\begin{document}
\maketitle {\bf Abstract:} In this paper, we present a new
extension of the famous Serrin's lower semicontinuity theorem for
the variational functional $\int_{\Omega}f(x,u,u')dx$,we prove its
lower semicontinuity in $W_{loc}^{1,1}(\Omega)$ with respect to
the strong $L_{loc}^{1}$ topology assuming that the integrand
$f(x,s,\xi)$ has the usual continuity on all the three variables
and the convexity property on the variable $\xi$ and the local
absolute continuity on the  variable $x$.\\

{\bf Keyword: } Lower semicontinuity, Serrin's theorem, strong convergence in $L^{1}$,
convex function, local absolute continuity.\\
{\bf 2002 Mathematics Subject Classification:} Primary 49J45,
Secondary 52A41.
\section{Introduction and Main Results}
\setcounter{equation}{0}

The aim of this paper is to give some new sufficient conditions
for lower semicontinuity with respect to the strong convergence in
$L_{loc}^{1}$ for functionals of integral type
\begin{equation}
F(u,\Omega)=\int_{\Omega}f(x,u(x),Du(x))dx,\label{(1.1)}
\end{equation}
where $\Omega$ is an open set of $R^{n}$, $u$ is in the Sobolev
Space$^{[1]}$ $W_{loc}^{1,1}(\Omega)=\{u:u\in L^{1}(K),Du\in
L^{1}(K), \forall K\subset\subset\Omega\}$, $Du$ denotes the
generalized gradient of $u$, and the integrand $f(x,s,\xi):\Omega
\times R \times R^{n}\rightarrow
[0,\infty)$ satisfies the following conditions:\\
(H1) $f$ is continuous in $\Omega \times R \times R^{n}$ and
$f(x,s,\xi)$ is convex in $\xi\in R^{n}$ for all $(x,s)\in \Omega
\times R$.\\
\indent The integral functional $F$ is called lower semicontinuous
in $W_{loc}^{1,1}(\Omega)$ with respect to the strong convergence
in $L_{loc}^{1}$, if for every $u_{m},u \in W_{loc}^{1,1}(\Omega)$
such that $u_{m}\rightarrow u $ in $L_{loc}^{1}$
 (where $u_{m}\rightarrow u $ in $L_{loc}^{1}$ means
 $\|u_{m}-u\|_{L^{1}(K)}\rightarrow 0 \ as \ m\rightarrow +\infty, \forall K\subset\subset\Omega$), then
\begin{equation}
\liminf_{m\rightarrow+\infty}F(u_{m},\Omega)\geq
F(u,\Omega).\label{(1.2)}
\end{equation}
It is well known that condition (H1) alone is not sufficient for
strong lower semicontinuity of the integral $F$ in (1.1) (see book
[12]). In addition to (H1),Serrin published in 1961 an
article$^{[13]}$ proposing some sufficient conditions for strong
lower semicontinuity. One of the most known
and celebrated Serrin's theorem on this subject is the following one.\\
{\bf Theorem 1.1$^{[13]}$ } Let f satisfy, in addition to (H1), one of
the
following conditions:\\
(a) $f(x, s, \xi)\rightarrow +\infty$ when $|\xi|\rightarrow
+\infty$, for all $(x,s)\in \Omega\times R $;\\
(b) $f(x, s, \xi)$ is strictly convex in $\xi \in R^{n}$
for all $(x,s)\in \Omega\times R $;\\
(c) the derivatives $f_{x}(x,s,\xi)$, $ f_{\xi}(x,s,\xi)$ and
$f_{\xi
x}(x,s,\xi)$ exist and are continuous.\\
Then $F(u,\Omega)$ is lower semicontinuous in
$W_{loc}^{1,1}(\Omega)$ with respect to the strong convergence in
$L_{loc}^{1}$.\\
\indent The conditions (a), (b) and (c) quoted above are clearly
independent, in the sense that we can find a continuous function
$f$ satisfying just one of them, but none of the others . However,
the proof of Theorem 1.1 is essentially the same for every
condition quoted above; indeed, the proof is based on an
approximation theorem for convex functions depending continuously
on parameters that can be applied, in particular, when f satisfies
one of conditions (a), (b) and (c). This fact suggests the
possibility to find a suitable condition weaker than one of
conditions (a), (b) and (c). Many attempts have been made to
weaken the assumptions on the integrand $f$, such as L. Ambrosio
in paper [2], V. De Cicco in his paper [3] and I. Fonseca in his
papers [6] and [7] proposed several generalizations of Theorem
1.1. In the papers [9] and
[10],  Gori prove the following theorems:\\
{\bf Theorem 1.2$^{[9]}$} Let us assume that $f$ satisfies (H1) and
also assume that, for every compact set $K \subset \Omega \times R
\times R^{n}$, there exists a constant $L = L(K)$ such that
\begin{equation}
|f_{\xi}(x_{1},s,\xi)-f_{\xi}(x_{2},s,\xi)|\leq L|x_{1}-x_{2}|, \
\forall (x_{1},s,\xi),(x_{2},s,\xi)\in K,\label{(1.3)}
\end{equation}
and, for every compact set $K_{1} \subset \Omega \times R$, there
exists a constant $L_{1} = L_{1}(K_{1})$ such that\\
\begin{equation}
|f_{\xi}(x,s,\xi)|\leq L_{1}, \ \forall (x,s)\in K_{1}, \ \forall
\xi \in R^{n},\label{(1.4)}
\end{equation}
\begin{equation}
|f_{\xi}(x,s,\xi_{1})-f_{\xi}(x,s,\xi_{2})|\leq
L_{1}|\xi_{1}-\xi_{2}|, \ \forall (x,s)\in K_{1}, \ \forall \xi_{1},
\xi_{2}\in R^{n}.\label{(1.5)}
\end{equation}
Then $F(u,\Omega)$ is lower semicontinuous in
$W_{loc}^{1,1}(\Omega)$ with respect to the strong convergence in
$L_{loc}^{1}$.\\
{\bf Theorem 1.3$^{[9]}$} Let f  satisfy (H1) and such that, for every
open set $\Omega'\times H \times K \subset \subset \Omega\times R
\times R^{n}$, there exists a constant $L=L_{\Omega'\times H \times
K}$ such that, for every $x_{1},x_{2}
\in \Omega'$, $s\in H$ and $\xi \in K$ ,\\
\begin{equation}
|f(x_{1},s,\xi)-f(x_{2},s,\xi)|\leq L|x_{1}-x_{2}|.\label{(1.6)}
\end{equation}
Then the functional $F(u,\Omega)$ is lower semicontinuous on
$W_{loc}^{1,1}(\Omega)$ with
respect to the $L_{loc}^{1}$ convergence.\\
\indent Condition (1.6) means that $f$ is locally Lipschitz continuous
with respect to $x$, locally with respect to $(s,\xi)$ and not
necessarily globally, that is, the Lipschitz constant is not uniform
for $(s,\xi)\in R \times R^{n}$. This is an improvement of (c) of
Serrin's Theorem 1.1  since, when only the gradient $\nabla_{x}f$
exists and is continuous, this implies the Lipschitz continuity of
$f$ with respect to $x$ on the
compact subsets of $\Omega\times R \times R^{n}$.\\
\indent Then a question arises that whether there are weaker
enough conditions more than locally Lipschitz continuous
condition? In this paper, we consider absolutely continuous
condition. Obviously, absolute continuity is weaker than Lipschitz
continuity. Is the local absolute continuity condition on the
integrand  enough for the lower semicontinuity of the variational
functional? The
following theorems  give a confirmed answer.\\
{\bf Theorem 1.4 } Let $\Omega \subset R$ be an open set, $f(x,s,\xi):\Omega\times R\times R\longrightarrow[0,+\infty)$
satisfy the following conditions:\\
(H1)    $f(x,s,\xi) $ is continuous on $\Omega\times R\times R$,  $f(x,s,\xi)$ is convex in $\xi \in R$ for all $(x,s)\in \Omega\times R$;\\
(H2)  $f_{\xi}(x,s,\xi) $ is continuous on $\Omega\times R\times R$, and for every compact set of $\Omega\times R\times R$,
$f_{\xi}(x,s,\xi) $ is absolutely continuous about $x$;\\
(H3) for every compact set $K_{1}\subseteq\Omega\times R$, there exists a constant $L_{1}=L_{1}(K_{1})$, such that\\
\begin{equation}
|f_{\xi}|\leq L_{1}, \ \forall(x,s)\in K_{1},\ \ \forall \xi\in
R,\label{(1.7)}
\end{equation}
\begin{equation}
|f_{\xi}(x,s,\xi_{1})-f_{\xi}(x,s,\xi_{2})|\leq
L_{1}|\xi_{1}-\xi_{2}|, \ \forall(x,s)\in K_{1},\ \ \forall
\xi_{1},\xi_{2}\in R.\label{(1.8)}
\end{equation}
Then the functional $F(u,\Omega)=\int_{\Omega}f(x,u(x),u'(x))dx$ is lower semicontinuous on
$W_{loc}^{1,1}(\Omega)$ with
respect to the strong convergence in $L_{loc}^{1}(\Omega)$.\\
{\bf Theorem 1.5 } Let $\Omega \subset R$ be an open set, $f(x,s,\xi):\Omega\times R\times R\longrightarrow[0,+\infty)$
satisfy (H1) and the following conditions:\\
(H4)  for every compact set $\Omega'\times H \times K\subseteq \Omega\times R \times R$,
$f(x,s,\xi) $ is absolutely continuous about $x$;\\
Then the functional $F(u,\Omega)$ is lower semicontinuous on
$W_{loc}^{1,1}(\Omega)$ with
respect to the strong convergence in $L_{loc}^{1}(\Omega)$.\\

\section{Some Lemmas}
\setcounter{equation}{0}

In this section, we collect some preliminary definitions and lemmas (see papers [1,4,8,11,14]) which will be
used in the sequel.\\
{\bf Definition 2.1}\ \ Let $\Omega\subset R^n$ be on open set, we
denote\\
 $$L^p_{loc}(\Omega )=\{u: \Omega\rightarrow R|\ \int_{\Omega^{\prime}}|u|^pdx<+\infty
,\forall \Omega^{\prime}\subset\subset\Omega\},$$
and\\
$$W^{1,1}_{loc}(\Omega )=\{u|\ u\in L^1_{loc}(\Omega ), Du\in
L^1_{loc}(\Omega )\}$$ {\bf Remark 2.1}\ \ Notice that
$u_m\rightarrow u$ in $L^1_{loc}(\Omega )$
implies that $u_m$ converges to $u$ in measure.\\
{\bf Definition 2.2}\ \ Let $f: [a,b]\rightarrow R$ be a real
function, if $\forall\varepsilon >0, \exists\delta >0$, such that
for any finite disjoint open interval $\{(a_i,b_i)\}^n_{i=1}$ on
$[a,b]$, when $\sum\limits^n_{i=1}(b_i-a_i)<\delta$, we have
$$\sum\limits^n_{i=1}|f(b_i)-f(a_i)|<\varepsilon ,$$
then we call $f(x)$ is a absolutely continuous function on $[a,b]$.\\
{\bf Remark 2.2}\ \ If $f(x)$ is Lipschitz continuous on $[a,b]$,
then $f(x)$ is absolutely continuous on $[a,b]$.\\
{\bf Lemma 2.1}\ \ Let $f(x)$ is a absolutely continuous function on
$[a,b]$, then $f(x)$ is almost everywhere differentiable on $[a,b]$
and $f^{\prime}(x)$ is a integrable function on $[a,b]$.\\
{\bf Lemma 2.2}(Lebesgue Dominated convergence theorem)\ \ Let
$(X,\mathcal{R},\mu)$ be a measure space, $f$ and $\{f_n\}(n\geq
1)$
be measurable functions on $E\in\mathcal{R}$, if\\
(1) $\{f_n\}$ is convergence in measure to $f$ on $E$;\\
(2) there exists a integrable function $h(x)$ on $E$, such that
$$|f_n(x)|\leq h(x),\ \ \ \ a.e. x\in E,$$
then $f(x), f_{n}(x)(n\geq 1)$ are integrable on $E$, and
$$\lim\limits_{n\rightarrow +\infty}\int_Ef_n(x)d\mu =\int_Ef(x)d\mu
.$$
{\bf Lemma 2.3}\ \ Let $f(x)$ be a measurable function on $E$, the
$f(x)$ is Lebesgue integrable on $E$ if and only if $|f(x)|$ is
Lebesgue integrable on $E$, and
$$|\int_Ef(x)dx|\leq\int_E|f(x)|dx$$
{\bf Lemma 2.4}\ \ Let $f(x)$ be Lebesgue integrable on $E, then
\forall\varepsilon >0, \exists\delta >0$, for every measurable
subset $A$ of $E$, when $m(A)<\delta$, we have
$$|\int_Af(x)dx|\leq\int_A|f(x)|dx<\varepsilon .$$
{\bf Definition 2.3}\ \ Let function $\eta (x): R^n\rightarrow
[0,1]$ satisfies $\eta\in C^{\infty}_c(R^n), {\rm supp}(\eta
)\subseteq B(0,1), \eta (-x)=\eta (x)$ and $\int_{R^n}\eta (x)dx=1$.
Giving a function $v: \Omega\rightarrow R$ and
$\varepsilon>0$, we define the convolution of $v$ with
step $\varepsilon$ as
$$v_{\varepsilon}(x)=\eta_{\varepsilon}*v(x)=\int_{R^n}\eta_{\varepsilon}(x-y)v(y)dy=\int_{R^n}\eta_{\varepsilon}(y)v(x-y)dy,$$
where
$\eta_{\varepsilon}(x)=\varepsilon^{-n}\eta\left(\frac{x}{\varepsilon}\right),
{\rm supp}(\eta )=\overline{\{x\in R^n,|\eta (x)\neq 0\}}.$\\
\indent  We have the following
approximations about the convolution of $v$:\\
{\bf Lemma 2.5$^{[1]}$}\ \ Let $\Omega\subset R^n$ be an open set, and a
function $v: \Omega\rightarrow R$, then \\
(i) If $v\in L^1_{loc}(\Omega )$ and ${\rm
supp}(v)\subset\subset\Omega $, then
$v_{\varepsilon}\in C^{\infty}_0(\Omega )$ provided $\varepsilon<dist({\rm
supp}(v),\partial\Omega)$, and
$$v_{\varepsilon}\rightarrow v\ {\rm in}\ L^1_{loc}(\Omega )\ {\rm
as}\ \varepsilon\rightarrow 0^{+}.$$
(ii) If $v\in L^p(\Omega )$ where $1\leq p<+\infty $, then
$v_{\varepsilon}\in L^{p}(\Omega )$. Also
$$v_{\varepsilon}\rightarrow v\ {\rm in}\ L^p_{loc}(\Omega)\ {\rm
as}\ \varepsilon\rightarrow 0^{+}.$$
\indent A function $f: R^n\rightarrow R$ is called convex if for
every $x,y\in R^n, \lambda\in (0,1)$, we have
$$f(\lambda x+(1-\lambda
)y)\leq\lambda f(x)+(1-\lambda )f(y)$$.

Now, we give some properties about convex functions:\\
{\bf Lemma 2.6$^{[4]}$.}\ \ Let $f_i: \Omega\rightarrow R,
\{f_i\}_{i\in N}$ be a sequence of functions ,\\
(i) if $f_i$ is convex, then $f=\sup_{i\in N}f_i$ is also convex;\\
(ii) if $f_i$ is lower semicontinuous, then $f=\sup_{i\in N}f_i$ is
also lower semicontinuous.
\indent The following approximation result was proved by De Giorgi$^{[8]}$.\\
{\bf Lemma 2.7$^{[8]}$}\ \ Let $U\subseteq R^d$ be an open set and
$f: U\times R^n\rightarrow [0,+\infty )$ be a continuous function
with compact support on $U$, such that, for every $t\in U,
f(t,\cdot )$ is convex on $R^n$. Then there exists a sequence
$\{\alpha_{q}\}^{\infty}_{q=1}\subseteq C^{\infty}_c(R^n),
\alpha_{q}\geq 0, \int_{R^n}\alpha_{q}dx=1, {\rm
supp}(\alpha_{q})\subseteq B(0,1)$.
 Let
$$a_{q}(t)=\int_{R^n}f(t,\xi )\{(n+1)\alpha_{q}(\xi )+<\nabla \alpha_{q}(\xi ),\xi
>\}d\xi ,$$
and
$$b_q(t)=-\int_{R^n}f(t,\xi )\nabla \alpha_q(\xi )d\xi .$$
Then
$$f_j(t,\xi )=\max\limits_{1\leq q\leq j}\{0,a_q(t)+<b_q(t),\xi
>\},j\in N,$$
satisfy the following results:\\
(i) for every $j\in N, f_j: U\times R^n\rightarrow [0,+\infty )$ is
a continuous function with compact support on $U$ such that, $\forall t\in U, f_j(t,\cdot )$ is convex on $R^n$.
Moreover, $\forall(t,\xi )\in U\times R^n, f_j(t,\xi )\leq f_{j+1}(t,\xi )$ and
$$f(t,\xi )=\sup\limits_{j\in N}f_j(t,\xi )$$\\
(ii) for every $j\in N$, there exists a constant $M_j>0$, such that,
$\forall (t,\xi )\in U\times R^n$,
$$|f_j(t,\xi )|\leq M_j(1+|\xi |),$$
and,  $\forall t\in U, \forall\xi_1,\xi_2\in R^n.$
$$|f_j(t,\xi_1)-f_j(t,\xi_2)|\leq M_j|\xi_1-\xi_2|.$$

\section{Proof of Theorem 1.4}
\setcounter{equation}{0}
We will divide into four steps to complete the proof of  theorem 1.4.\\
{\bf Step 1:}\ \ Let $\{\beta_i(x,s)\}_{i\in N}$ be a sequence of
smooth functions satisfying \\
(1) there exists a compact set $\Omega^{\prime}\times
H\subset\subset\Omega\times R$, such that $\beta_i(x,s)=0, \forall
(x,s)\in (\Omega\backslash\Omega^{\prime})\times (R\backslash H),$\\
(2) for every $i\in N, \beta_i(x,s)\leq\beta_{i+1}(x,s),\forall
(x,s)\in\Omega'\times H$,\\
(3) $\lim\limits_{i\rightarrow +\infty}\beta_i(x,s)=1,\forall
(x,s)\in\Omega'\times H$.\\
\indent
Let
$$f_i(x,s,\xi )=\beta_i(x,s)f(x,s,\xi ),\ \ \ \ \ i=1,2,\cdots .$$
\indent It is clear that, for each $i\in N, f_i$  satisfies all
the hypothesis in theorem 1.4 and also vanishes if $(x,s)$ is
outside $\Omega'\times H$. Thus
$$\lim\limits_{i\rightarrow +\infty}f_i(x,s,\xi )=f(x,s,\xi ),\ \ \
\ \forall (x,s,\xi )\in\Omega'\times H\times R,$$ and
$$f_i(x,s,\xi )\leq f_{i+1}(x,s,\xi )\leq f(x,s,\xi ),\forall i\in
N, \forall (x,s,\xi )\in\Omega'\times H\times R.$$ By Levi lemma,
we have
$$\lim\limits_{i\rightarrow
+\infty}\int_{\Omega'}f_i(x,s,\xi
)dx=\int_{\Omega'}f(x,s,\xi )d\xi .$$
Thus, without loss of
generality, we can assume that there exists an open set
$\Omega'\times H\subset\subset\Omega\times R$, such that
\begin{equation}
f(x,s,\xi )=0,\ \ \ \ \forall (x,s,\xi )\in (\Omega
\backslash\Omega')\times (R\backslash H)\times R.\label{(3.1)}
\end{equation}
\indent Let $u_m, u\in W^{1,1}_{loc}(\Omega )$ such that $u_m\rightarrow u$
in $L^1_{loc}(\Omega )$. We will prove that
$$\liminf\limits_{m\rightarrow +\infty} F(u_m,\Omega )\geq
F(u,\Omega ).$$
 Without loss of generality, we can assume that
$$\liminf\limits_{m\rightarrow +\infty} F(u_m,\Omega
)=\lim\limits_{m\rightarrow +\infty}F(u_m,\Omega )<+\infty .$$
By (3.1), we have $F(u_m,\Omega )=F(u_m,\Omega'), F(u,\Omega
)=F(u,\Omega')$, thus we will only prove the following inequality:
\begin{equation}
\lim\limits_{m\rightarrow +\infty}F(u_m,\Omega')\geq
F(u,\Omega').\label{(3.2)}
\end{equation}
{\bf Step 2:}\ \  Let $\eta_{\varepsilon}\in C^{\infty}_c(R)$ be a
mollifier and for $\epsilon>0$,define
\begin{equation}
v_{\varepsilon}(x)=\eta_{\varepsilon}*u(x)=\int_{\Omega}\eta_{\varepsilon}(x-y)u(y)dy,\
\ \ \ x\in [\Omega_{\varepsilon}],\label{(3.3)}
\end{equation}
where $\left[\Omega_{\varepsilon}\right]\triangleq \{x\in\Omega
:{\rm dist}(x,\partial\Omega )>\varepsilon\}$. We have
\begin{eqnarray}
[u_{\varepsilon}(x)]'&=&[\eta_{\varepsilon}*u(x)]_x=[\int_{\Omega}\eta_{\varepsilon}(x-y)u(y)dy]_x\nonumber\\
&=&\int_{\Omega}[\eta_{\varepsilon}(x-y)]_xu(y)dy=\int_{\Omega}-[\eta_{\varepsilon}(x-y)]_yu(y)dy\label{(3.4)}\\
&=&\int_{B(x,\varepsilon
)}\eta_{\varepsilon}(x-y)[u(y)]_ydy=[u']_{\varepsilon}(x),\
\ \ \ x\in\Omega_{\varepsilon}.\nonumber
\end{eqnarray}
In the following, we denote the derivative of $u_{\varepsilon}$ as
$u'_{\varepsilon}$. When $u\in W^{1,1}_{loc}(\Omega )$, we know
$u'\in L^1_{loc}(\Omega )$. By Lemma 2.5, we know
$u'_{\varepsilon}\in C^{\infty}_0(\Omega )$ and
\begin{equation}
u'_{\varepsilon}\rightarrow u'\ {\rm in}\ L^1_{loc}(\Omega )\ {\rm
as}\ \varepsilon\rightarrow 0^{+},\label{(3.5)}
\end{equation}
i.e., $\forall\delta >0, \ \exists\epsilon>0$, such that
\begin{equation}
\int_{\Omega'}|u'_{\varepsilon}-u'|dx<\delta .\label{(3.6)}
\end{equation}
\indent New we estimate the  difference for the integrand values
on different points:
\begin{eqnarray}
f(x,u_m,u'_m)-f(x,u,u')&=&f(x,u_m,u'_m)-f(x,u_m,u'_{\varepsilon})\nonumber\\
&+&f(x,u_m,u'_{\varepsilon})-f(x,u,u'_{\varepsilon})\label{(3.7)}\\
&+&f(x,u,u'_{\varepsilon})-f(x,u,u').\nonumber
\end{eqnarray}
By the convexity of $f(x,s,\xi )$ with respect to $\xi$, we have
\begin{equation}
f(x,u_m,u'_m)-f(x,u_m,u'_{\varepsilon})\geq
f_{\xi}(x,u_m,u'_{\varepsilon})\cdot(u'_m-u'_{\varepsilon}).\label{(3.8)}
\end{equation}
By (\ref{(3.8)}), we have
\begin{eqnarray}
f(x,u_m,u'_m)-f(x,u_m,u'_{\varepsilon})&\geq & f_{\xi}(x,u_m,u'_{\varepsilon})\cdot
u'_m-f_{\xi}(x,u_m,u'_{\varepsilon})\cdot u'_{\varepsilon}\nonumber\\
&=&f_{\xi}(x,u_m,u'_{\varepsilon})\cdot
u'_m-f_{\xi}(x,u,u'_{\varepsilon})\cdot u'\nonumber\\
&+&f_{\xi}(x,u,u'_{\varepsilon})\cdot
u'-f_{\xi}(x,u_m,u'_{\varepsilon})\cdot u'_{\varepsilon}\label{(3.9)}\\
&=&f_{\xi}(x,u_m,u'_{\varepsilon})\cdot
u'_m-f_{\xi}(x,u,u'_{\varepsilon})\cdot u'\nonumber\\
&+& f_{\xi}(x,u,u'_{\varepsilon})\cdot
(u'-u'_{\varepsilon})\nonumber\\
&+& [f_{\xi}(x,u,u'_{\varepsilon})-f_{\xi}(x,u_m,u'_{\varepsilon})]\cdot
u'_{\varepsilon}.\nonumber
\end{eqnarray}
By (\ref{(3.7)}) and (\ref{(3.9)}), we have
\begin{eqnarray}
f(x,u_m,u'_m)-f(x,u,u')&\geq & f_{\xi}(x,u_m,u'_{\varepsilon}) \cdot
u'_m-f_{\xi}(x,u,u'_{\varepsilon})\cdot u'\nonumber\\
& +&f_{\xi}(x,u,u'_{\varepsilon})\cdot
(u'-u'_{\varepsilon})\nonumber\\
&+& [f_{\xi}(x,u,u'_{\varepsilon})-f_{\xi}(x,u_m,u'_{\varepsilon})]\cdot
u'_{\varepsilon}\label{(3.10)}\\
&+&f(x,u_m,u'_{\varepsilon})-f(x,u,u'_{\varepsilon})\nonumber\\
&+& f(x,u,u'_{\varepsilon})-f(x,u,u').\nonumber
\end{eqnarray}
This implies
\begin{eqnarray}
&&\int_{\Omega'}[f(x,u_m,u'_m)-f(x,u,u')]dx\nonumber\\
&\geq &\int_{\Omega'}[f_{\xi}(x,u_m,u'_{\varepsilon})\cdot
u'_m-f_{\xi}(x,u,u'_{\xi})\cdot
u']dx\nonumber\\
&+&\int_{\Omega'}[f_{\xi}(x,u,u'_{\varepsilon})\cdot
(u'-u'_{\varepsilon})]dx \nonumber\\
&+&\int_{\Omega'}[f_{\xi}(x,u,u'_{\varepsilon})-f_{\xi}(x,u_m,u'_{\varepsilon})]\cdot
u'_{\varepsilon}dx\label{(3.11)}\\
&+&\int_{\Omega'}[f(x,u_m,u'_{\varepsilon})-f(x,u,u'_{\varepsilon})]dx\nonumber\\
&+&\int_{\Omega'}[f(x,u,u'_{\varepsilon})-f(x,u,u')]dx.\nonumber
\end{eqnarray}
{\bf Step 3:}\ \ Now, we estimate the right side of inequality
(\ref{(3.11)}).\\
\indent  By (\ref{(1.7)}) and (\ref{(3.6)}), we have
\begin{equation}
\int_{\Omega'}[f_{\xi}(x,u,u'_{\varepsilon})\cdot
(u'-u'_{\varepsilon})]dx\geq
-L_1\int_{\Omega'}|u'-u'_{\varepsilon}|dx\geq
-L_1\delta .\label{(3.12)}
\end{equation}
Thus
\begin{equation}
\lim\limits_{\varepsilon\rightarrow
0}\int_{\Omega'}[f_{\xi}(x,u,u'_{\varepsilon})\cdot
(u'-u'_{\varepsilon})]dx\geq 0\label{(3.13)}
\end{equation}
Since $f(x,s,\xi )$ and $f_{\xi}(x,s,\xi )$ are continuous functions
, they are bounded functions on compact subset. By remark 2.1 and Lebesgue
dominated convergence theorem, we obtain
\begin{equation}
\lim\limits_{m\rightarrow
+\infty}\int_{\Omega'}[f_{\xi}(x,u,u'_{\varepsilon})-f_{\xi}(x,u_m,u'_{\varepsilon})]\cdot
u'_{\varepsilon}dx=0,\label{(3.14)}
\end{equation}
and
\begin{equation}
\lim\limits_{m\rightarrow
+\infty}\int_{\Omega'}[f(x,u_m,u'_{\varepsilon})-f(x,u,u'_{\varepsilon})]dx=0.\label{(3.15)}
\end{equation}
Now, we will prove
\begin{equation}
\lim\limits_{\varepsilon\rightarrow
0}\int_{\Omega'}[f(x,u,u'_{\varepsilon})-f(x,u,u')]dx\geq
0\label{(3.16)}
\end{equation}
By lemma 2.7, there exists a sequence of non-negative continuous
functions $f_j(x,s,\xi )\ (j\in N)$, such that $f_j(x,s,\xi )$ is
convex on $\xi$, and $\forall (x,s,\xi )\in\Omega'\times
H\times R$,
\begin{eqnarray}
&&f_j(x,s,\xi )\leq f_{j+1}(x,s,\xi ),\label{(3.17)}\\
&&f(x,s,\xi )=\sup\limits_{j\in N}f_j(x,s,\xi ),\label{(3.18)}\\
&&|f_j(x,s,\xi_1)-f_j(x,s,\xi_2)|\leq
M_j|\xi_1-\xi_2|.\label{(3.19)}
\end{eqnarray}
By Levi lemma, we obtain
\begin{equation}
\lim\limits_{j\rightarrow
+\infty}\int_{\Omega'}f_j(x,u,u'_{\varepsilon})dx=\int_{\Omega'}f(x,u,u'_{\varepsilon})dx,\label{(3.20)}
\end{equation}
and
\begin{equation}
\lim\limits_{j\rightarrow
+\infty}\int_{\Omega'}f_j(x,u,u')dx=\int_{\Omega'}f(x,u,u')dx.\label{(3.21)}
\end{equation}
In order to prove (\ref{(3.16)}), we only need to prove
\begin{equation}
\lim\limits_{\varepsilon\rightarrow
0}\int_{\Omega'}[f_j(x,u,u'_{\varepsilon})-f_j(x,u,u')]dx\geq
0,\ \forall j\in N.\label{(3.22)}
\end{equation}
By (\ref{(3.20)}) and (\ref{(3.21)}), we have
$$\int_{\Omega'}[f_j(x,u,u'_{\varepsilon})-f_j(x,u,u')]dx\geq
-M_j\int_{\Omega'}|u'_{\varepsilon}-u'|dx\geq
-M_j \delta .$$ Thus (\ref{(3.16)}) holds.\\
{\bf Step 4:} Now, we need to prove
\begin{equation}
\lim\limits_{m\rightarrow
+\infty}\int_{\Omega'}[f_{\xi}(x,u_m,u'_{\varepsilon})\cdot
u'_m-f_{\xi}(x,u,u'_{\varepsilon})\cdot
u']dx=0. \label{(3.23)}
\end{equation}

Let
\begin{eqnarray}
&&g(x,s)\triangleq
f_{\xi}(x,s,u'_{\varepsilon}),\label{(3.24)}\\
&&G_m(x)\triangleq\int^{u_m(x)}_{u(x)}g(x,s)ds.\label{(3.25)}
\end{eqnarray}
By condition $(H_2)$, $ g(x,s)$ is a absolutely continuous function on
$x$. By Lemma 2.1, $g(x,s)$ is almost everywhere differentiable,
i.e. $\frac{\partial g}{\partial x}$ exists almost everywhere.
Derivating the both sides of (\ref{(3.25)}), we obtain
\begin{equation}
G^{\prime}_m(x)=g(x,u_m)\cdot u'_m-g(x,u)\cdot
u'+\int^{u_m(x)}_{u(x)}\frac{\partial g}{\partial x}dx,\ \ \
\ a.e. x\in\Omega'.\label{(3.26)}
\end{equation}
Because $G_m(x)$ vanishes outside $\Omega'$, we obtain
\begin{equation}
\int_{\Omega'}G^{\prime}_m(x)dx=0. \label{(3.27)}
\end{equation}
By (\ref{(3.26)}), we have
\begin{eqnarray}
&&|\int_{\Omega'}[f_{\xi}(x,u_m,u'_{\varepsilon})\cdot
u'_m-f_{\xi}(x,u,u'_{\varepsilon})\cdot
u']dx|\nonumber\\
&=&\left|\int_{\Omega'}[g(x,u_m)\cdot
u'_m-g(x,u)\cdot
u']dx\right|\label{(3.28)}\\
&=&|-\int_{\Omega'}\int^{u_m(x)}_{u(x)}\frac{\partial
g}{\partial x}dsdx|\leq\int_{D_m}\left|\frac{\partial g}{\partial
x}\right|dxds, \nonumber
\end{eqnarray}
where
$$D_m=\{(x,s)\in\Omega'\times H|\ \min\{u_m(x),u(x)\}\leq
s(x)\leq\max\{u_m(x),u(x)\}\}.$$ We note
\begin{equation}
|D_m|=|\int_{\Omega'}\int^{u_m}_udsdx|\leq\int_{\Omega'}|u_m-u|dx\rightarrow
0\ (m\rightarrow +\infty ).\label{(3.29)}
\end{equation}
By Fubini theorem, we have
\begin{equation}
\int_{\Omega\times R}|\frac{\partial g}{\partial
x}|dxds=\int_Hds\int_{\Omega'}\left|\frac{\partial
g}{\partial x}\right|dx\label{(3.30)}
\end{equation}
Since $g(x,s)$ is absolutely continuous about $x$, $\frac{\partial
g}{\partial x}$ is integrable and absolutely integrable, i.e.
\begin{equation}
\int_{\Omega'}\left|\frac{\partial g}{\partial
x}\right|dx<+\infty .\label{(3.31)}
\end{equation}
Thus
\begin{equation}
\int_{\Omega\times R}\left|\frac{\partial g}{\partial
x}\right|dxds<+\infty\label{(3.32)}
\end{equation}
Because of the absolute continuity of integral, we have
\begin{equation}
\lim\limits_{m\rightarrow +\infty}\int_{D_m}\left|\frac{\partial
g}{\partial x}\right|dxds=0\label{(3.33)}
\end{equation}
By (\ref{(3.28)}), we obtain
\begin{equation}
\lim\limits_{m\rightarrow
+\infty}|\int_{\Omega'}[f_{\xi}(x,u_m,u'_{\varepsilon})\cdot
u'_m-f_{\xi}(x,u,u'_{\varepsilon})\cdot
u']dx|=0\label{(3.34)}
\end{equation}
Thus we have proved (\ref{(3.23)}). By (\ref{(3.13)})-(\ref{(3.16)})
and (\ref{(3.23)}), we have
\begin{equation}
\lim\limits_{m\rightarrow
+\infty}\int_{\Omega'}[f(x,u_m,u'_m)-f(x,u,u')]dx\geq
0\label{(3.35)}
\end{equation}
Thus we deduce that the functional $F(u,\Omega )$ is lower
semicontinuous on $W^{1,1}_{loc}(\Omega )$ with respect to the
strong convergence in $L^1_{loc}(\Omega )$, which does
complete the proof.\\
\section{Proof of Theorem 1.5}
\setcounter{equation}{0}

In order to prove Theorem 1.5, we will verify all the conditions
in Theorem 1.4 from the assumptions in Theorem 1.5.
Now we will
divide into four steps to complete the
proof of Theorem 1.5:\\
{\bf Step 1:}\ \ Similar to the first step of the proof in Theorem
1.4, without loss of generality, we assume that the integrand
$f(x,s,\xi)$ vanishes outside a compact subset of $\Omega\times
R$. Thus we assume that there exists an open set $\Omega'\times
H\subset\subset\Omega\times R$, such that
\begin{equation}
f(x,s,\xi )\equiv 0,\ \ \ \ \forall (x,s,\xi )\in
(\Omega\backslash\Omega')\times (R\backslash H)\times
R.\label{(4.1)}
\end{equation}
\indent Let $u_m,u\in W^{1,1}_{loc}(\Omega )$, such that $u_m\rightarrow u$
in $L^1_{loc}(\Omega )$, we need to prove
\begin{equation}
\lim\limits_{m\rightarrow +\infty}F(u_m,\Omega')\geq
F(u,\Omega').\label{(4.2)}
\end{equation}
By lemma 2.7, there exists a functions sequence $\{f_j(x,s,\xi
)\}_{j\in N}$, such that $\forall j\in N$, $f_j$ is a continuous
function on $\Omega'\times H\subset\subset\Omega\times R$;
$\forall (x,s)\in\Omega'\times H$, $f_j(x,s,\cdot )$ is convex on
$R$, and $\forall (x,s,\xi )\in\Omega'\times H\times R$,
\begin{eqnarray}
&&f_j(x,s,\xi )\leq f_{j+1}(x,s,\xi )\label{(4.3)}\\
&&f(x,s,\xi )=\sup\limits_{j\in N}f_j(x,s,\xi )\label{(4.4)}\\
&&|f_j(x,s,\xi_1)-f_j(x,s,\xi_2)|\leq
M_j|\xi_1-\xi_2|,(x,s)\in\Omega'\times H, \xi_1,\xi_2\in
R.\label{(4.5)}
\end{eqnarray}

Let $\eta_{\varepsilon}\in C^{\infty}_c(R)(0<\varepsilon <<1)$ be a
mollifier and define the $f_{j,\varepsilon}=f_j*\eta_{\varepsilon}$,
i.e.
\begin{equation}
f_{j,\varepsilon}(x,s,\xi )=\int_Rf_j(x,s,\xi
-z)\eta_{\varepsilon}(z)dz.\label{(4.6)}
\end{equation}
By (\ref{(4.5)}), we have
\begin{eqnarray}
&&|f_{j,\varepsilon}(x,s,\xi )-f_j(x,s,\xi )|\nonumber\\
&=&|\int_Rf_j(x,s,\xi -z)\eta_{\varepsilon}(z)dz-\int_Rf_j(x,s,\xi
)\eta_{\varepsilon}(z)dz|\nonumber\\
&\leq&\int_R|f_j(x,s,\xi -z)-f_j(x,s,\xi
)|\eta_{\varepsilon}(z)dz\label{(4.7)}\\
&\leq&\int_{{\rm
supp}\eta_{\varepsilon}}M_j|z|\cdot\eta_{\varepsilon}(z)dz\leq
M_j\cdot\varepsilon. \nonumber
\end{eqnarray}
Choosing $\varepsilon =\varepsilon_j=\frac{1}{jM_j}\rightarrow 0$. By
(\ref{(4.7)}), we have \begin{equation} |f_{j,\varepsilon_j}(x,s,\xi
)-f_j(x,s,\xi )|\leq M_j\varepsilon_j=\frac{1}{j}.\label{(4.8)}
\end{equation}
So
\begin{equation}
f_j(x,s,\xi )-\frac{2}{j}\leq f_{j,\varepsilon_j}(x,s,\xi
)-\frac{1}{j}\leq f_j(x,s,\xi )\leq f(x,s,\xi )\label{(4.9)}
\end{equation}
By (\ref{(4.3)}), (\ref{(4.4)}) and Levi lemma, we have
\begin{equation}
\lim\limits_{i\rightarrow
+\infty}\int_{\Omega'}f_j(x,u(x),u'(x))dx=\int_{\Omega'}f(x,u(x),u'(x))dx.\label{(4.10)}
\end{equation}
Denote
\begin{equation}
F_j(u,\Omega')=\int_{\Omega'}[f_{j,\varepsilon_j}(x,u(x),u'(x))-\frac{1}{j}]dx\label{(4.11)}
\end{equation}
By (\ref{(4.9)})-(\ref{(4.11)}), we have
\begin{equation}
\lim\limits_{j\rightarrow +\infty}F_j(u,\Omega')=F(u,\Omega'
)=\int_{\Omega'}f(x,u(x),u'(x))dx.\label{(4.12)}
\end{equation}
Obviously,
$$F_j(u,\Omega')\leq F(u,\Omega'),\ \ \ \ \forall
j\in N.$$ Thus
\begin{equation}
\sup\limits_{j\in
N}F_j(u,\Omega')=F(u,\Omega').\label{(4.13)}
\end{equation}
By Lemma 2.6, in order to prove that $F(u,\Omega')$ is lower
semicontinuous on $W^{1,1}_{loc}(\Omega )$ with respect to the
strong convergence in $L^{\prime}_{loc}(\Omega )$, it's enough
that we  prove lower semicontinuity with respect to the strong
convergence in $L^1_{loc}(\Omega )$ for the functional sequence
$\{F_j(u,\Omega')\}_{j\in N}$  on $W^{1,1}_{loc}(\Omega )$ .\\
{\bf Step 2:} \ In order to prove that $\forall j\in N,
F_j(u,\Omega')$ is lower semicontinuous on
$W^{1,1}_{loc}(\Omega )$ with respect to the strong convergence in
$L^1_{loc}(\Omega )$, we will prove that $\forall j\in N$, the
integrand $f_{j,\varepsilon_j}(x,u(x),u'(x))$ satisfy all
conditions of theorem 1.4.\\
\indent  $\forall\xi_1,\xi_2\in R, 0<\lambda <1$,
by the convexity of $f_j(x,s,\cdot )$ on $R$, we have
\begin{equation}
f_j(x,s,\lambda\xi_1+(1-\lambda )\xi_2)\leq\lambda
f_j(x,s,\xi_1)+(1-\lambda )f_j(x,s,\xi_2).\label{(4.14)}
\end{equation}
Thus
\begin{eqnarray}
&&f_{j,\varepsilon_j}(x,s,\lambda\xi_1+(1-\lambda )\xi_2)\nonumber\\
&=&\int_Rf_j(x,s,\lambda\xi_1+(1-\lambda
)\xi_2-z)\eta_{\varepsilon_j}(z)dz\nonumber\\
&=&\int_Rf_j(x,s,\lambda (\xi_1-z)+(1-\lambda
)(\xi_2-z))\cdot\eta_{\varepsilon_j}(z)dz\label{(4.15)}\\
&\leq&\lambda\int_Rf_j(x,s,(\xi_1-z))\cdot\eta_{\varepsilon_j}(z)dz+(1-\lambda
)\int_Rf_j(x,s,(\xi_2-z))\cdot\eta_{\varepsilon_j}(z)dz\nonumber\\
&=&\lambda f_{j,\varepsilon_j}(x,s,\xi_1)+(1-\lambda
)f_{j,\varepsilon_j}(x,s,\xi_2).\nonumber
\end{eqnarray}
Thus $f_{j,\varepsilon_j}$ satisfy $(H_1)$.\\
{\bf Step 3:}\ \ $\forall (x,s)\in\Omega'\times H, \xi_1,\xi_2\in
R$, by (\ref{(4.5)}), we have
\begin{eqnarray}
&&|f_{j,\varepsilon_j}(x,s,\xi_1)-f_{j,\varepsilon_j}(x,s,\xi_2)|\nonumber\\
&=&|\int_R[f_j(x,s,\xi_1-z)\eta_{\varepsilon_j}(z)-f_j(x,s,\xi_{2}-z)\eta_{\varepsilon_j}(z)]dz|\nonumber\\
&\leq&\int_R|f_j(x,s,\xi_1-z)-f_j(x,s,\xi_2-z)|\cdot\eta_{\varepsilon_j}(z)dz\label{(4.16)}\\
&\leq&\int_{{\rm
supp}\eta_{\varepsilon}}M_j|\xi_1-\xi_2|\eta_{\varepsilon_j}(z)dz\leq
M_j|\xi_1-\xi_2|.\nonumber
\end{eqnarray}
Thus
\begin{equation}
\left|\frac{\partial f_{j,\varepsilon_j}}{\partial\xi}\right|\leq
M_j.\label{(4.17)}
\end{equation}
So $f_{j,\varepsilon_j}$ satisfies (1.7) in the condition (H3) of Theorem 1.4.\\
\indent  Now, we will prove $f_{j,\varepsilon_j}$ satisfies (1.8)
in the condition (H3) of Theorem 1.4. By ${\rm
supp}(\eta_{\varepsilon_j})\subseteq B(0,\varepsilon_j)$,we have

\begin{eqnarray}
\frac{\partial f_{j,\varepsilon_j}}{\partial\xi}(x,s,\xi
)&=&\int_R\frac{\partial f_j(x,s,\xi
-z)}{\partial\xi}\cdot\eta_{\varepsilon_j}(z)dz\nonumber\\
&=&-\int_R\frac{\partial f_j(x,s,\xi -z)}{\partial
z}\cdot\eta_{\varepsilon_j}(z)dz\label{(4.18)}\\
&=&\int_Rf_j(x,s,\xi
-z)\frac{\partial\eta_{\varepsilon_j}(z)}{\partial
z}dz.\nonumber
\end{eqnarray}

By (\ref{(4.5)}) and (\ref{(4.18)}), we have
\begin{eqnarray}
&&\left|\frac{\partial
f_{j,\varepsilon_j}}{\partial\xi}(x,s,\xi_1)-\frac{\partial
f_{j,\varepsilon_j}}{\partial\xi}(x,s,\xi_2)\right|\nonumber\\
&\leq&\int_R\left|f_j(x,s,\xi_1-z)-f_j(x,s,\xi_2-z)\right|\cdot\left|\frac{\partial\eta_{\varepsilon_j}(z)}{\partial
z}\right|dz\label{(4.19)}\\
&\leq&M_j|\xi_1-\xi_2|\cdot\int_R\left|\frac{\partial\eta_{\varepsilon_j}(z)}{\partial
z}\right|dz=L_jM_j|\xi_1-\xi_2|,\nonumber
\end{eqnarray}
where
\begin{equation}
L_j=\int_R\left|\frac{\partial\eta_{\varepsilon_j}(z)}{\partial
z}\right|dz\label{(4.20)}
\end{equation}
is a constant depending on $\varepsilon_j$. Thus
$f_{j,\varepsilon_j}$ satisfies (1.8). So we have proved that
$f_{j,\varepsilon_j}$ satisfies $(H_3)$.\\
{\bf Step 4:}\ \ Next we will prove that $f_{j,\varepsilon_j}$
satisfies condition $(H2)$.\\
\indent By the condition $(H4)$, for every compact subset
$\Omega'\times H\times K$, $f(x,s,\xi )$ is absolutely continuous
about $x$, i.e. $\forall\varepsilon_0>0,\ \exists\delta >0$ such
that for any finite disjoint open interval $\{(x_i,y_j)\}^n_{i=1}$
in $\Omega'$, when $\Sigma^n_{i=1}(y_i-x_i)<\delta$, we have
\begin{equation}
\sum\limits^n_{i=1}|f(y_i,s,\xi )-f(x_i,s,\xi
)|<\varepsilon_0\label{(4.21)}
\end{equation}
By Lemma 2.7, there exists continuous functions sequence
$\{f_j(x,s,\xi )\}_{i\in N}$,  $\forall j\in N$, $\forall
(x,s)\in\Omega'\times H, \ f_j(x,s,\cdot )$ is convex on $R$,
and $\forall (x,s,\xi )\in\Omega'\times H\times R$, we have
\begin{equation}
f_j(x,s,\xi )=\max\limits_{1\leq q\leq
j}\{0,a_q(x,s)+b_q(x,s)\xi\},\ j\in N.\label{(4.22)}
\end{equation}
where
\begin{eqnarray}
&&a_q(x,s)=\int_Rf(x,s,\xi )\left[2\eta_q(\xi
)+\xi\frac{\partial\eta_q(\xi )}{\partial\xi}\right]d\xi
,\label{(4.23)}\\
&&b_q(x,s)=-\int_Rf(x,s,\xi )\frac{\partial\eta_q(\xi
)}{\partial\xi}d\xi ,\label{(4.24)}
\end{eqnarray}
and $\eta_q\in C^{\infty}_c(R)\ (q\in N)$ are mollifiers
satisfying $\eta_q\geq 0, \ \int_R\eta_q(\xi )d\xi =1$ and ${\rm
supp}(\eta_q)\subseteq B(0,1)$, $\forall j\in N$. By
(\ref{(4.22)}), without of loss generality, we assume that there
exists $l\in \{1,\cdots,j\}$, such that
\begin{equation}
f_j(x,s,\xi )=a_l(x,s)+b_l(x,s)\cdot\xi\label{(4.25)}
\end{equation}
where $a_l,b_l$ are given by (\ref{(4.23)})-(\ref{(4.24)}). By
(\ref{(4.21)}), we obtain
 \begin{eqnarray}
&&\sum\limits^n_{i=1}|a_l(y_i,s)-a_l(x_i,s)|\nonumber\\
&=&\sum\limits^n_{i=1}|\int_R[f(y_i,s,\xi )-f(x_i,s,\xi )]\cdot
[2\eta_l(\xi )+\xi\frac{\partial\eta_l(\xi )}{\partial\xi}]d\xi
|\nonumber\\
&\leq&\int_R\sum\limits^n_{i=1}|f(y_i,s,\xi )-f(x_i,s,\xi )|\cdot
[2\eta_l(\xi )+|\xi\frac{\partial\eta_l(\xi
)}{\partial\xi}|]d\xi\label{(4.26)}\\
&\leq&\varepsilon_0\int_{B(0,1)}[2\eta_l(\xi
)+|\xi\frac{\partial\eta_l(\xi )}{\partial\xi}|]d\xi\leq
(2+A_l)\cdot\varepsilon_0,\nonumber
\end{eqnarray}
where $$A_l=\int_{B(0,1)}|\frac{\partial\eta_l(\xi
)}{\partial\xi}|d\xi$$ is a constant.
Similar to the above proof,
we have
\begin{eqnarray}
&&\sum\limits^n_{i=1}|b_l(y_i,s)-b_l(x_i,s)|\nonumber\\
&=&\sum\limits^n_{i=1}|\int_R[f(y_i,s,\xi )-f(x_i,s,\xi
)]\cdot\frac{\partial\eta_l(\xi )}{\partial\xi}d\xi
|\nonumber\\
&\leq&\int_R\sum\limits^n_{i=1}|f(y_i,s,\xi )-f(x_i,s,\xi
)|\cdot\left|\frac{\partial\eta_l(\xi
)}{\partial\xi}\right|d\xi\label{(4.27)}\\
&\leq&\varepsilon_0\int_{B(0,1)}\left|\frac{\partial\eta_l(\xi
)}{\partial\xi}\right|d\xi\leq A_l\cdot\varepsilon_0.\nonumber
\end{eqnarray}
Thus \begin{eqnarray} &&\sum\limits^n_{i=1}|f_j(y_i,s,\xi
)-f_j(x_i,s,\xi )|\nonumber\\
&=&\sum\limits^n_{i=1}|a_l(y_i,s)-a_l(x_i,s)+[b_l(y_i,s)-b_l(x_i,s)]\cdot\xi
|\nonumber\\
&\leq&\sum\limits^n_{i=1}|a_l(y_i,s)-a_l(x_i,s)|+\sum\limits^n_{i=1}|b_l(y_i,s)-b_l(x_i,s)|\cdot
|\xi |\label{(4.28)}\\
&\leq&(2+A_l)\varepsilon_0+A_l\varepsilon_0K_1=(2+A_l+A_lK_1)\varepsilon_0\triangleq\sigma.\nonumber
\end{eqnarray}
Since $\xi$ varies on a compact set, then
$K_1=\sup\limits_{\xi}\{|\xi |\}<+\infty$. Choosing
$\varepsilon_0$  sufficient small so that $\sigma$ is enough
small. Thus $f_j(x,s,\xi )$ is absolutely continuous about $x$ on
any compact subset of $\Omega\times R\times R$. By (\ref{(4.6)})
and (\ref{(4.28)}), we have \begin{eqnarray}
&&\sum\limits^n_{i=1}|f_{j,\varepsilon_j}(y_i,s,\xi
)-f_{j,\varepsilon_j}(x_i,s,\xi
)|\nonumber\\
&=&\sum\limits^n_{i=1}\left|\int_R[f_j(y_i,s,\xi
-z)\eta_{\varepsilon_j}(z)-f_j(x_i,s,\xi-
z)\eta_{\varepsilon_j}(z)]dz\right|\nonumber\\
&\leq&\int_R\sum\limits^n_{i=1}|f_j(y_i,s,\xi -z)-f_j(x_i,s,\xi
-z)|\cdot\eta_{\varepsilon_j}(z)dz\label{(4.29)}\\
&\leq&\sigma\cdot\int_{B(0,\varepsilon_j)}\eta_{\varepsilon_j}(z)dz=\sigma
.\nonumber \end{eqnarray}
By (\ref{(4.18)}) and (\ref{(4.29)}), we
obtain
 \begin{eqnarray}
  &&\sum\limits^n_{i=1}\left|\frac{\partial
f_{j,\varepsilon_j}}{\partial\xi}(y_i,s,\xi)-\frac{\partial
f_{j,\varepsilon_j}}{\partial\xi}(x_i,s,\xi)\right|\nonumber\\
&\leq&\int_R\sum\limits^n_{i=1}|f_j(y_i,s,\xi -z)-f_j(x_i,s,\xi
-z)|\cdot\left|\frac{\partial\eta_{\varepsilon_j}(z)}{\partial
z}\right|dz\label{(4.30)}\\
&\leq&\sigma\int_R\left|\frac{\partial\eta_{\varepsilon_j(z)}}{\partial
z}\right|dz=L_j\sigma\nonumber
\end{eqnarray}
where $L_j$ are constants depending on $\varepsilon_j$ and given by
(\ref{(4.20)}) ($\forall j\in N$). By (\ref{(4.30)}), for every compact subset on
$\Omega\times R\times R, \frac{\partial
f_{j,\varepsilon_j}}{\partial\xi}$ is absolutely continuous about $x$.
Thus $f_{j,\varepsilon_j}$ satisfy condition
$(H_2)$.

Now,we have proved $f_{j,\varepsilon_j}$ satisfies all conditions
in Theorem 1.4, so $F_j(u,\Omega')$ is lower semicontinuous in
$W^{1,1}_{loc}(\Omega )$ with respect to the strong convergence in
$L^1_{loc}(\Omega )$. Thus $F(u,\Omega )$ has the same lower
semicontinuity. This completes the proof of Theorem 1.5.\\


{}
\end{document}